\documentclass{article}
\usepackage{amssymb, amsfonts, amsmath, amsthm}
\usepackage[russian]{babel}
\newtheorem{theorem}{Theorem}

\newtheorem{lemma}{Lemma}
\newtheorem{cor}[theorem]{Corollary}

\newtheorem{rem}{Remark}
\DeclareMathOperator{\intr}{int}  % внутренность
\DeclareMathOperator{\dil}{dil}    % dilatation

\newcommand{\ep}{\epsilon}

\newcommand{\al}{\alpha}

\newcommand{\be}{\beta}

\newcommand{\tr}{\triangle}
\newcommand{\an}{\angle}
\newcommand{\D}{\partial}

\newcommand{\ov}{\overline}
\newcommand{\R}{\mathbb R}   % might be the standard \mathbf style

\begin{document}
\title{Bi-Lipschitz equivalent Alexandrov surfaces, I}
\author{A.~Belenkiy, Yu.~Burago
\footnote{The second author was partly supported by grants  RFBR 02-01-00090,
SS-1914.2003, and CRDF RM1-2381-ST-02}
}
%\address{  }
%\email{burago@pdmi.ras.ru}
%\date{}
\maketitle
\section{Basic definitions and statements}

Recently M.~Bonk and U.~Lang    \cite{BonkLang} proved, that if a
complete  Riemannian manifold $M$ is homeomorphic to the plane
$\R^2$ and satisfies the conditions  $\int_M K^+dS <2\pi$,\,
$\int_M K^-dS<\infty$, then Lipschitz distance $d_L(M,\R^2)$
between  $M$ and $\R^2$ satisfies the inequality

\centerline{$\displaystyle{d_L(M,\R^2)\leq\ln\Big(\frac{2\pi+\int_M K^-dS }
{2\pi-\int_M K^+dS}\Big)^{\frac{1}{2}}}.$}

 This inequality is sharp if curvature does not change its sign.
Here $K^+=\max\{K,0\}$, $K^-=\max\{-K,0\}$, and $dS$ is the area element. We will remind
 definition of Lipschitz metric a little later.

In fact, this result was obtained in  \cite{BonkLang} for the class of Alexandrov
surfaces  more wide than the class of Riemannian manifolds.

This paper is inspired by the paper  \cite{BonkLang} mentioned above; we will investigate
surfaces which are not necessary simply connected.
In contrast to the case of surfaces homeomorphic to $\R^2$, in more general cases
we do not have any standard models.
By this reason we try to estimate the Lipschitz distance
between two homeomorphic surfaces. Our estimate occurs to be far
from the optimal one (and in this part we even restrict ourselves
proving finiteness of   distances).
Our consideration is naturally divided  into two parts: asymptotic at infinity and
study of compact surfaces.

Our readers supposed to be familiar with the basic notions of two dimensional manifolds
of bounded total (integral) curvature theory. Its expositions can be
found, for instance,  in
\cite{AlexZalg} and \cite{Reshetnyak}.

Hereafter  Alexandrov surface means a {\em complete} two dimensional
manifold of bounded curvature with a boundary; the boundary (which may be empty)
is supposed to consist of finite number of curves with finite variation of turn.

We  introduce the following notations:
let  $M$ be an Alexandrov surface with a metric  $d$,\;
$\omega$ be its curvature, which is a signed measure,
 $\omega^+$, $\,\omega^-$ be the positive and negative parts of the curvature, and
  $\Omega=\omega^++\omega^-$ be the variation of the curvature.
 For any Riemannian manifold,
$\omega^+=\int_M K^+dS$,\; $\omega^-=\int_M K^-dS$.

A point $p$ carrying curvature $2\pi$ and a boundary point carrying  turn $\pi$
are called peak points.

Recall that dilatation  $\dil f$ of a Lipschitz map
$f\colon X\to Y$, where $(X,\,d_X)$ and $(Y,\,d_Y)$ are metric spaces,
is defined by the equality
$$
  \dil f = \sup_{x,x'\in X,\,x\not= y} \frac{d_Y(f(x),f(x'))}{d_X(x,x')}.
$$
A homeomorphism  $f$ is called bi-Lipschitz if both maps,
 $f$ and $f^{-1}$, are Lipschitz ones.
The impression
$$
  d_L(X,Y) = \inf_{f\colon X\to Y} \ln\bigl(\max \{\dil(f),\,
\dil(f^{-1})\}\bigr)
$$
\noindent is called Lipschitz distance between  $X$ and $Y$; here infimum is taken
over all Lipschitz homeomorphisms  $f\colon X\to~Y$.
Metric spaces $X$, $Y$ are bi-Lipschitz equivalent if and only if
$d_L(X,Y)<\infty$.

\begin{theorem}
\label{theorem-1} Let two compact Alexandrov surfaces are homeomorphic one to another,
have no peak points. Then these surfaces are bi-Lipschitz equivalent one to another.
\end{theorem}

 This theorem is trivial for
Riemannian manifolds because two homeomorphic smooth 2-manifolds are always
diffeomorphic. Theorem \ref{theorem-1} does not give an upper estimate
for Lipschitz distance via some finite set of geometric characteristics of the surfaces,
like their diameters, total curvatures, systolic constants, etc.
However such an estimate does exist; we are going to publish
this result in a separate paper.
\medskip

Let $T$ be an end; i.e., an Alexandrov  surface homeomorphic to a closed disk with
the removed center   and such that for every sequence of
points  $p_i\in T$ with images in the disk converging to the center, the
condition  $d(a, p_i)\to\infty$ holds for
$i\to\infty$; here
$a$ is any fixed point. Such a sequence is called diverging
(or going to infinity) one.

Let us denote the turn of the end boundary $\D T$ by $\sigma$.
We call  the value $v(T)=-\omega (T)-\sigma$  a growth speed of the end  $T$.
The Cohn-Vossen inequality says that this value is not negative.
The growth speed of an end is positive if and only if
$\displaystyle{\lim_{i\to\infty}\frac{l(\gamma_i)}{d(a, p_i)}>0}$,
where $l(\gamma_i)$ is the length of the shortest noncontractible loop
with the  vertex $p_i$.
Under condition  $\Omega (T)<\infty$, this limit is well-defined and is not greater than 2.

As it was proved in
(\cite{Huber},\, \cite{Verner}),
a complete Alexandrov surface satisfying the condition
 $\omega^- (M)<\infty$ is homeomorphic to a closed
 surface with  finite number of removed points. Ends are appropriate closed
 neighborhoods of these points. So, any Alexandrov surface can be cut into
 a compact part and some ends.

\begin{theorem}
\label{theorem-2} Let two complete Alexandrov surfaces
$M_1$, $M_2$ be homeomorphic one to another, satisfy the condition
$\omega^-(M_i)<\infty$ and contain neither  points with curvature
$2\pi$, nor boundary points with turn
$\pi$. Also, let  all ends of these surfaces have nonzero growth speed.
Then these surfaces are bi-Lipschitz equ\-iva\-lent.
\end{theorem}

\begin{rem} {\rm  We call two ends to be equivalent if their intersection
contains an end. It is not difficult to see that equivalent ends have equal growth speeds.
This means that the theorem above does not depend on how are cut surfaces into
compact parts and ends.
}
\end{rem}

Simple examples show that ends having zero speed (even having no peak points)
are not necessary bi-Lipschitz
equivalent. For instance, none of the surfaces obtained
by rotation of the
following graphs are  bi-Lipschitz equivalent:

 $\{y=\sqrt{x},\,x\geq 1\}$,\qquad $\{y=1, \, x\geq 0$\},\qquad $\{y=e^{-x}, \, x\geq 0\}$.

\noindent Note that surfaces of rotation obtained from the last graph and the graph of
the function
$\{y=e^{-2x}, \, x\geq 0\}$ are bi-Lipschitz equivalent.

Also note that an end having zero growth speed can not be bi-Lipschitz
equivalent to an end with nonzero growth speed.

An end is rotationally symmetric if
its isometry group contains a subgroup whose restriction to the boundary is
transitive.
The problem of  bi-Lipschitz equivalence of ends can be  reduced to the similar
problem for rotationally symmetric ends.   More precisely, the following theorem
takes place.

\begin{theorem}
\label{theorem-3}
Every end $T$ satisfying the condition $\omega^- (T)<\infty$ and
having no peak points is bi-Lipschitz equivalent to a rotationally
symmetric end with  the same growth speed.
\end{theorem}

\begin{rem}
\label{Lang-Bonk ends} {\rm If  $v(T)>0$, then the end $T$ contains a smaller
end $T_1$, such that the Lipschitz distance between $T_1$ and a plane with a disk
of length $l=l(\partial T_1)$ removed can be estimated from above by a value
depending on $v$ и $l$ only.
}
\end{rem}

This remark can be proved by a minor modification of arguments from \cite{BonkLang};
so we omit the proof.

\bigskip

We see that the problem of  bi-Lipschitz classification of ends with
zero growth speed is reduced to a problem for functions in one variable.
Intuitively, one can imagine a zero speed end as ``having a peak point
(with curvature $2\pi$) at infinity''. From this point of view, it is naturally to expect
 analogy between classification of zero speed ends and bi-Lipschitz
classification of neighborhoods of finite peak points.
 Note that for smooth surfaces with isolated singularities the
 classification problem for neighborhoods of peak points was considered by
 D.~Grieser \cite{Grieser}.

\medskip

{\bf Sketch of the proof.}
%%%Временное условное название???

The idea of our proof of Theorem \ref{theorem-1} is simple. Here is the sketch of the
proof.
Preliminary   a simply connected region  $\tr$ bounded by
a simple closed curve with three selected points was called a generalized triangle.
These points are  vertices of the triangle, intervals of the curve between
vertices being its sides.
We always assume that lengths
of the sides satisfy the triangle inequality, the sides are geodesic broken lines.
Actually we will consider only generalized triangles with small variations of
curvature and  turn of its sides. We will often omit the word ``generalized''.

A partition of a 2-manifold $M$ is a set of generalized
triangles in $M$, whose interiors do not overlap and whose union is the whole $M$.
A triangulation of $M$ is a partition such that for any two triangles of it
their intersection is either  a side or a vertex.

We will partition both  Alexandrov  surfaces   $M_1$, $M_2$  into
generalized triangles in such a way that these triangles  have  small variations of
curvature and  turn of its sides, but angles are separated from zero.
In particular, we include all points carrying essential portion of curvature
in the set of vertices of the triangles . We prove that each triangle of the partition
is bi-Lipschitz equivalent to its comparison triangle; i.e., to a planar
triangles with the same side lengths; one can choose such a  bi-Lipschitz mapping
to preserve length of sides. For this we use the method introduced by
I.~Bakelman \cite{Bakelman} for constructing Tchebysev coordinates in Alexandrov
 surfaces.
 Replacing each generalized triangle of our partition by its comparison  triangle
and attaching last planar triangles together in according with the same combinatorial
scheme, we will get two  surfaces $P_1,\; P_2$ equipped with polyhedral metrics.
In according to our construction, each polyhedron  $P_i$ is
bi-Lipschitz equivalent to the surface  $M_i$, $i=1,2$. Finally, note that the
polyhedra  $P_1$ and $P_2$ are bi-Lipschitz equivalent one to another.
This ends the proof.

Theorem  \ref{theorem-2} follows from Theorem  \ref{theorem-1}  and the lemma below.

\begin{lemma}
\label{BL-lemma}
Every end  $T$ with nonzero growth speed, without peak points  and with
$\omega^-(T)<\infty$ is  bi-Lipschitz equivalent to $\R^2$ with a disk removed
(as always, we suppose that boundary of the end is a curve
with finite variation of turn and no peak points on the end and its boundary).
Moreover, there is a  bi-Lipschitz equivalence inducing an affine map on the boundary
of the end.
\end{lemma}

In case when curvature and turn of the boundary are not great, the lemma can be
proved by a minor modification of arguments from  \cite{BonkLang}; to do this
it is sufficient to suppose that  $\omega^+(T)+\tau^+(\partial T)<\pi$, where
$\tau^+(\partial T)$ is the  positive  turn of the end boundary (from the side of
the end).

To prove the lemma in the general case, it is sufficient to cut the end into an
annular and an end satisfying the last condition and to apply
Theorem  \ref{theorem-1} to the annular.

Lemma  \ref{BL-lemma}  also follows from Theorem \ref{theorem-3} and Theorem
\ref{theorem-1}. Indeed,  Theorem \ref{theorem-3} easily implies that
every end with nonzero growth speed is bi-Lipschitz
equivalent to a cone over a  circle of appropriate length with a round neighborhood
of the vertex removed. The last surface is obviously bi-Lipschitz
equivalent to the plane with a disk removed.

To prove Theorem  \ref{theorem-3}, we will use again a modification
of I.~Bakelman's construction \cite{Bakelman}
.

\section {Triangles having small curvature.}
%Slightly curved triangles --- variant????
We begin with simple statements about planar triangles.
Let  $\triangle ABC$ be a flat triangle. Denote by   $a,\;b,\;c$  its sides
opposite to the angles  $\angle A$, $\angle B$, $\angle C$,
correspondingly. We will denote side lengths by the same letters as sides.

\begin{lemma}
\label{affin-1}
Let planar triangles   $\triangle ABC$ and  $\triangle\ov A\ov B\ov C$ be
 such that $b=\ov b$, $c=\ov c$, $\angle A\leq\angle \ov A\leq L\cdot \angle A$,
 and $L\angle \ov A-\angle A\leq\pi(L-1)$, where  $L>1$.
Then the optimal  bi-Lipschitz constant of  the affine transformation
mapping  $\triangle ABC$ to
$\triangle \ov A \ov B \ov C$ is not greater than ~$L$.
\end{lemma}

It is sufficient  to find eigenvalues of the  affine transformation under consi\-de\-ration
to prove the lemma.

\rightline{$\blacksquare$}
\medskip

\begin{cor}
\label{affin-2}
Let  triangles  $\triangle ABC$ and   $\triangle \ov A \ov B \ov C$
satisfy the conditions:
$$
\begin{array}{cc}
  L^{-1}\le c/\ov c\le L, &  L^{-1}\le b/\ov b\le L, \\
  \ep\le\an A\le\pi -\ep, & \ep\le\an \ov A\le\pi -\ep.
\end{array}
$$
Then the optimal bi-Lipschitz constant of the affine transformation mapping
$\tr ABC$ on $\tr \ov A\ov B\ov C$ is not greater than
 $\displaystyle{(L\frac{\pi}{2\epsilon})^2}$.
\end{cor}

\begin{lemma}
\label{affin-3}
Let  triangles  $\triangle ABC$ and   $\triangle \ov A \ov B \ov C$
satisfy the condition:
 $b=\ov b$,  $c=\ov c$,\; $D\in BC$. Denote by  $\ov D$ the image of the  point $D$
 under affine transformation mapping  $\tr ABC$ onto  $\tr \ov A \ov B \ov C$.

Let
$\an BAD=\al,\;\an CAD=\be,\;\an \ov B\ov A\ov D=\al_1,\;\an \ov C\ov A\ov D=\be_1$.
Suppose that
$\displaystyle{\al <\frac{\pi}{2}, \; \be < \frac{\pi}{2}}$ and
$$
0<\angle \ov B \ov A \ov C- \al<\frac{\pi}{2}, \;
0<\angle \ov B \ov A \ov C- \be<\frac{\pi}{2}.
$$

Then
$$
|\al-\al_1| \le |\angle BAC - \angle \ov B \ov A\ov C|,\;
 |\be-\be_1| \le |\angle BAC - \angle \ov B \ov A \ov C|.
$$
 \end{lemma}
\medskip

{\bf Proof.}
Suppose  that $\angle \ov B \ov A \ov C>\angle BAC$. Let us show that in this case
 $\al_1 \ge \al, \; \be_1 \ge \be $.
After these equalities
  $\al + \be = \angle BAC,\;
   \al_1 + \be_1 = \angle \ov B \ov A \ov C $,
  will imply the lemma.

Since $\ov D$ is the image of   $D$,

 $$
 \frac{\sin\al}{\sin\be}=\frac{\sin\al_1}{\sin\be_1}.
 $$
Suppose  that $\al_1 < \al$. Consider the point
 $Z$ in  $\ov B\ov C$  such that
$\angle \ov B \ov A Z = \al$.
Then
$$
\frac{\sin\al}{\sin\be}<\frac{\sin(\ov B \ov A Z)}{\sin(\ov C \ov A Z)},
$$
therefore $\sin\be>\sin(\angle {\ov B \ov A \ov C}-\al )$. But
$0<\be<(\angle \ov B \ov A \ov C- \al)< \frac{\pi}{2}$.
This means that $\sin\be<\sin(\angle {\ov B \ov A \ov C}- \al)$.
Contradiction.

The case $\angle \ov B \ov A \ov C<\angle BAC$ can be considered analogously.

\rightline{$\blacksquare$}
\medskip

Now let us go back to Alexandrov surfaces.
For a generalized  triangle $\tr\subset M$, we denote

\centerline{$\Omega'(\tr)=\Omega(\intr(\tr))+\sigma(\tr)$,}

\noindent where $\sigma(\tr)$ is the sum  of variations of  turn of the triangle sides
 from inside.

\begin{lemma}
\label{small-tr}

For every  $\theta>0$, there exist numbers  $\delta(\theta) >0$ and $L(\theta)>1$
with the following property:
if every angle of a generalized triangle $\tr ABC$ is not less than
$\theta$, and if   $\Omega'(\tr ABC)<\delta(\theta)$, then there exists a
bi-Lipschitz map of the triangle $ABC$ on its comparison  triangle with the constant
$L(\theta)$, whose restriction on the boundary of the  triangle preserves lengths.
\end{lemma}

This lemma looks like to be obvious,  and the fact we could find neither an
appropriate reference nor a very simple proof a little surprising seems to be
very surprising.
\medskip

 {\bf Proof.} 1. It is sufficient to prove our lemma only for polyhedral metrics.
So we suppose that our metric is a  polyhedral one. We will use a construction
which is a minor modification of I.~Bakelman's one, see \cite{Bakelman}.
Suppose that our   $\tr ABC$ (equipped with a polyhedral metric) is a part
of a complete surface  $M$ homeomorphic to $\R^2$ and  flat outside
the triangle.
 (For that, let us cut out a comparison triangle $\tr \ov A \ov B\ov C$ from $\R^2$
and then attach  $\tr ABC$ instead of $\tr \ov A \ov B\ov C$.)
Then we extend sides  $ \ov A \ov B$ and $ \ov A \ov C$ as rays
$\ov AB^*$ and $\ov AC^*$. Denote the sector  $B^*\ov AC^*$ by $S$.

We will partition the sector $S$ (more precisely,  some region of it
containing  $\tr ABC$) onto flat parallelograms coming into
contact one with  another
along the whole sides in such a way  that four  parallelograms adjust
to each vertex (except ones belonging the boundary of $S$).
This allows immediately to
introduce Tchebysev coordinates in the sector $S$. These coordinates give
(after some \linebreak additional deformation straightening the side $BC$) required
bi-Lipschitz map.

2. To simplify exposition, suppose that  $\tr ABC$ is an ordinary
(not \linebreak generalized) triangle with zero turn of its sides. The general case
differs from this model   by nonessential details only.
 Note that the difference of corresponding angles of
triangles $\tr ABC$ and $\tr \ov A \ov B\ov C$ is not greater than some function
of  $\Omega'(\tr ABC)$ (this function can be given explicitly),
which goes to zero \linebreak together with  $\Omega'$. The
proof of that is standard and based
on the theorem about ``arc and chord''
%span ???
(\cite{AlexZalg}, Lemma 5 of Chapter 9),
 %%о дугe и хорде,
the Gauss--Bonnet theorem, and comparison theorem for (nongeneralized) triangles.

As a result, we can choose   $\delta$ to be so small (smallness depends on
$\theta$ only), that every angle of the
triangle    $\tr \ov A \ov B\ov C$ is not less than
$\frac{2}{3}\theta$.

3. Let   $\an \ov C$ be the smallest angle of the triangle, and
$\an \ov A$ be the  greatest one.
Then angles $\an \ov B$ and  $\an \ov C$ are acute and each of them do not exceed
$\frac12 \pi-\frac{1}{3}\theta$.

4. Let us consider the line  $l$ containing the
ray  $AB^*$ and begin to  shift it (continuously) parallel to itself
inside  the sector  $S$. Initially this line will cut
 a flat  ``oblique'' semi-strip from $S$. We will continue this process
% this movement??
until vertices of the metric appear in  $l$ for the first time.
Denote by  $A_{11}$ the intersection of $l$ and $AC^*$, and let
$ A_{12}, \dots , A_{1n_1}$ be vertices of the metric which appears
in $l\cap S$, numbered ``from left to right''.

From every point $A_{1i}$ we emulate
% eject, emit  ???
a geodesic $A_{1i}B_{1i}$ outside the cut off semi-strip and so
that  $|A_{1i}B_{1i}|=|A_{1j}B_{1j}|$ for all $i, \,j$ and besides
$\an B_{1i}A_{1i}A_{1,i+1}+\an A_{1i} A_{1,i+1}B_{1,i+1}=\pi$.
If the length of $|A_{1i}B_{1i}|$
is sufficiently small, then we get a strip consisting of flat parallelograms
$ A_{1i}B_{1i}B_{1,i+1}A_{1,i+1}$. Let us choose $|A_{1i}B_{1i}|$
so that the vertexes of the
metric appear for the first time on the broken line $B_{11}B_{12} B_{1,3}\dots$.
Following this process we will get the partition
of the sector $S$ by parallelograms.
(the last one in this set of parallelograms is supposed to be an infinite oblique semi-strip)
To make parallelograms to be adjacent along the whole sides, it is sufficient to slit
some of them into more narrow
parallelograms.

5. If $\delta<\theta$ the process described above can be continued until we exhaust
all sector $S$  and  it  will produce bijective
mapping of the sector $S$ onto the first quadrant $S_0$ of the plane with the
oblique coordinates $Ouv$, where the coordinate angle $\an uOv$ equals to the angle
$\an A=\alpha$ of the sector $S$.

Indeed,  there could be only one obstacle: namely, at some step the broken line
$\Lambda_k=A_{k1}, A_{k2}\dots$ might  touch itself or the ray $AB^*$. In both cases
we would obtain a closed broken line
bounding a region $G$ homeomorphic to some disk. Applying the Gauss--Bonnet formula to $G$,
we would come to a contradiction with smallness of the curvature of our triangle
(in compare with its angles).

6. Now we can introduce Tchebychev coordinates in $S$. For the coordinates $(u, v)$
 of a point $p$ we'll take the lengths of the broken lines connecting this point with the
 rays $AC^*$,  $AB^*$, such that in every parallelogram crossed by these lines the latter
 ones go along the intervals parallel the corresponding sides of the parallelogram.
 Let $S_0$ be the first coordinate quadrant of the plane $\R^2$, and
 $\varphi\colon S\to S_0$ be the coordinate mapping constructed according to this plan.

7. On $S_0$ we define  the  metric by   the   linear element

$$
ds^2=du^2 + 2\cos\tau(u,v)dudv + dv^2,
$$
where
$$
\tau(u,v)=\alpha - \omega(\varphi^{-1}(D_{uv}))
    $$
and $D_{uv}$ is the parallelogram $[0\le u'<u,\, 0 \le v'<v]$.

This linear element makes  $S_0$ to be a metric space
$(S_0, d_1)$.

8. The map $\varphi$ is an isometry of $S$ to $(S_0, d_1)$. It is
not difficult to check this consequently along all parallelograms
of our partition of $S$.

Consider, in addition, the standard flat metric $d_2$ defined by
the linear element $ds^2=du^2 + 2\cos \alpha dudv + dv^2$ on
$S_0$. The map $\;\mbox{id}\colon (S_0,d_1)\to (S_0,d_2)$ is linear
on every parallelogram. The inequality
$|\tau(x,y)-\alpha|\le\Omega'(S)<\delta$ shows that this map is
bi-Lipschitz  with constant $(\delta \mu)$. Here and further we
denote by $\mu$ positive constants (may be different) which depend
on   $\theta $ only.

The image of  the side $BC$ of the triangle $\tr ABC$ under the
map $\varphi_1=\mbox{id}\circ\varphi$ is a broken line
$\Lambda=Y_0Y_1\dots Y_l$, where $Y_0=B'=\varphi_1B$,
$Y_l=C'=\varphi_1C$.    We want to prove, that there is a
$\mu$-bi-Lipschitz transformation $\zeta$, which move the region
$Q$, bounded by this broken and the shortests  $A'B'$, $A'C'$ (where
$A'=\varphi_1 A$ is the vertex of the sector $S_0$) onto the  triangle $\tr
A'B'C'$. After that, it will be  sufficient to apply Lemma
\ref{affin-1} to the flat triangles  $\tr A'B'C'$ and
$\tr \ov A\ov B\ov C$ to finish the proof of the lemma.

The transformation $\zeta$ is defined as  follows. Denote by $X_i$ the
intersection of the ray $A'Y_i$  and the side $B'C'$ of the triangle $\tr A'B'C'$.
Now we map affine each triangle $\tr A'Y_iY_{i+1}$ onto the corresponding
 triangle $\tr A'X_iX_{i+1}$.  Let us show that it gives a
 $\mu$-bi-Lipschitz map we need.

Indeed, every interval $Y_iY_{i+1}$ of the broken $\Lambda$ is an
affine image of an interval located in one of the parallelograms.
By our construction, parallelograms in $(S_0,d_2)$ are disposed by
``horizontal'' rows. Adding ``horizontal'' and ``vertical'' rays we can
consider each point  $Y_i$ to be a vertex of a parallelogram. As $\delta$
is small in compare with $\theta$, every interval $Y_iY_{i+1}$ is a diagonal
of the corresponding parallelogram.
Let a ray $N_i$ be the ``upper'' bound of
$k$-th row, so that  $Y_i\in N_i$, and
 $\tilde N_i= \varphi_1^{-1}(N_i)$. It is not
difficult to see that variation of turn of broken line $\tilde N_i$
is not greater than $\Omega'(M)<\delta$ (this is a rough estimate).
 Let us apply the Gauss--Bonnet formula to the region bounded by the shortest
 $AB$, intervals of the shortests $AC,\;BC$ and the curve $\tilde N_i$. This
 gives immediately that  the angle between  $\tilde N_i$ and
 the shortest $BC$ is different from the angle$\an B$ not greater than on $2\delta$.
 That means that the angle between $\tilde N_i$ and starting at $\tilde N_i$ edge of
 the broken line $\Lambda$, and so also the angle between the same edge and the ray
 $A'B'$, is differ from the angle $\an B'$ of the  triangle $\tr A'B'C'$ not greater
 than on  $3\delta$ (that follows easily from Lemma  \ref{affin-3}).
Now it is not difficult to calculate that

$$
|\frac{A'Y_i}{A'X_i}-1|<\frac{10\delta}{\alpha}
$$

\noindent (under the condition, that $\delta$  is small in comparison with
$\alpha$). Applying  Lemma \ref{affin-2}  to the couples
$\tr A'X_iX_{i+1}$ and $\tr A'Y_iY_{i+1}$, we get that
triangle $\tr A'B'C'$ and $\tr ABC$ are bi-Lipschitz equivalent with
some constant  $\mu$.Unfortunately  the restriction of the map we constructed
is not an isometry for the side $BC$ . Nevertheless this restriction changes distances
not greater than in
$\mu$ times; therefore it is possible to correct our map in every  triangle $\tr
A'X_iX_{i+1}$ in  such a way that it  remains to be a bi-Lipschitz equivalence
(with some constant  $\mu '$) but  becomes to be an isometry on the boundary of the
triangle. So we got the required map.

\rightline{$\blacksquare$}

\begin{lemma}
\label{small-1}
There exists $\xi>0$ such that if a generalized triangle $\tr ABC$ in  Alexandrov
space  satisfies the conditions:
$$
\angle ABC \ge \frac{\pi}{5},\; \angle ACB \ge  \frac{\pi}{5},\;
\an BAC>0,\; \Omega'(\tr ABC)<\xi,
$$
then there exists a  bi-Lipschitz map of the  triangle $\tr ABC$  onto its
comparison triangle whose restriction to the boundary maps
vertices to the corresponding
vertices and preserves lengths.
\end{lemma}

In contrast to the previous lemma, now we allow one of angles to be arbitrary small;
but from another side now  we can not estimate bi-Lipschitz constant.

 {\bf Proof.}
To  prove our Theorem, we cut
our  triangle into
triangles in  such a way, that each of them is bi-Lipschitz equivalent to
its comparison triangle.
We suppose that $\alpha =\an BAC<\frac{1}{100}$, otherwise Lemma \ref{small-tr} would say that
our triangle is  $L(\frac{1}{100})$-bi-Lipschitz equivalent to its comparison triangle
 $\tr\ov A\ov B\ov C$. Now we choose
 $\chi=\min\{\delta (\frac{1}{100}\pi ), \, \frac{1}{1000}\}$, see
Lemma \ref{small-tr}.

It is easy to see that there exist points
$A_1, C_1, B_1$ on the sides $BC,\; AC$   and $AB$, correspondingly,
such that

\centerline{$|AB_1|=|A C_1|,\; |B A_1|=|BB_1|,\; |CA_1|=|C C_1|.$}

\noindent Consider shortests (w.r.t.  the induced
metric of the triangle) connecting these points. They partition $\tr ABC$
onto four triangles
(see an explanation below). Let us repeat the same construction for the
triangle $\tr AB_1C_1$; i.e., choose points $A_2, C_2, B_2$ in the sides
 $B_1C_1,\; AC_1$     и $AB_1$, such that

\centerline{$|AB_2|=|A C_2|,\; |B_1 A_2|=|B_1B_2|,\; |C_1A_2|=|C_1 C_2|.$}

Let us continue this process. It is not difficult to calculate that all the angles of all
the triangles we obtained are separated from zero and $\pi$, for instance, they lie
between  $0,05\pi$ and  $0,95\pi$.
Therefore, according to Lemma  \ref{small-tr}, all such  triangles, except the triangle
$\tr AB_iC_i$, are $L(0,05\pi)$-bi-Lipschitz equivalent to their comparison
triangles.

Note that the partition process of the triangle can be continued up to infinity  and
$\sum |B_iC_i|=\sum (|B_iB_{i+1}|+|C_iC_{i+1}|)\le |AB|+|AC|$,
so that $|B_iC_i|\to 0$, when  $i\to\infty$, and $B_i,\,C_i\to A$, as $\alpha=\an BAC>0$.
It follows that there exists an $i=i_0$ such that
$\Omega'(\tr AB_{i_0}C_{i_0})<\delta (\alpha)$ and hence
$\tr AB_{i_0}C_{i_0}$  is bi-Lipschitz equivalent to its comparison triangle
(again according to Lemma \ref{small-tr}).

We complete the proof of the theorem using the backward induction on $i$:
$\tr AB_{i_0}C_{i_0}$ is bi-Lipschitz equivalent to its comparison triangle.
Suppose, that the same is true for $\tr AB_{i+1}C_{i+1}$, where $i\le i_0$, and prove for
$\tr AB_iC_i$. The latter is cut into 4 triangles each being  bi-Lipschitz equivalent
to its comparison triangle. Consider the analogous partition  for the comparison triangle
$\tr \ov A\ov B_i\ov C_i$. As the curvature  of  $\tr ABC$ is small being compared with
the angles of the triangles under consideration (except $\an BAC$), triangles of the
partition of the triangle $\tr \ov A\ov B_i\ov C_i$ are almost the same as comparison
triangles for the corres\-pond\-ing triangles of our partition
of  $\tr  A B_i C_i$. That is why it is easy
to construct bi-Lipschitz mappings for these couples of triangles. With this the
 proof is completed.

\rightline{$\blacksquare$}

\section{Compact  Alexandrov surfaces}

To prove Theorem \ref{theorem-1} we need the following

\begin{lemma}
\label{triangulation}
For any  $\xi>0$ there exists a partition by generalize triangles
of every
compact Alexandrov surface $M$ without peak points such that
for any triangle $\tr$ of this partition  has positive angles and
$\Omega'(\tr)<\xi$.
\end{lemma}

{\bf Proof.} First of all, we   triangulate a small neighborhood
of each point $p$ such that $\Omega (p)\leq\frac12 \xi$, where $ \xi$ was
chosen according to lemma \ref{small-1}.  We choose such neighborhoods
to be bounded by broken lines and  not overlapping. It
is easy to triangulate these neighborhoods to satisfy conditions of the lemma.

To  partition  the remain part of the surface  $M$ we use
Theorems 2  from Chapter 3 of the book \cite{AlexZalg}:
Any compact subset of $M$ with a polyhedral boundary can be covered by a set of
arbitrary small pairwise nonoverlapping simple triangles such that in every triangle
there is no side equals to the sum of two other sides.

It is clear that these triangles can be chosen to be so small that for each of them
$\Omega'<\frac12 \xi$. All that is left, is to deform these triangles to turn them into
generalized triangles to remove the zero angles and not to break the other suppositions
of the Lemma.
It is sufficient to remove one zero angle and to use  induction over the number of
zero angles. Suppose that there is at least one zero angle. As we have no peak points,
we can find two adjacent angles,
say, $\an BOA$ and $\an AOC$, the first one being equal to zero and the second
one being nonzero. Note, that the curve  $AOC$ may  either consist of two sides or
be a whole side of a generalized triangle.
Now it is sufficient to replace a very short initial interval of
the curve $OA$ by a two-component broken line  $ODE$, where $E\in OA$, lying
in the sector $\an AOC$, by a very small but nonzero angle with $OA$ and such that
variation of the turn of the broken line $ODEA$  exceeds
variation of the turn of the shortest line $OA$ very slightly.

\rightline{$\blacksquare$}

{\bf Proof of Theorem \ref{theorem-1}}.

Let us consider  a partition of a compact Alexandrov surface $M$ satisfying  Lemma
\ref{triangulation}. We want to prove  that every triangle of this partition
is bi-Lipschitz equivalent to
a flat triangle with the same side lengths. Besides, corresponding  bi-Lipschitz maps
can be chosen so that their restrictions  to the boundaries of the
triangles have to be  isometries.
For a triangle with two angles not less than  $ \frac15\pi$ each, this is true by Lemma
\ref{small-1}. Therefore we can assume that the triangle $\tr$ under consideration has
two angles, each of which is less than $\frac15\pi$. Then the third angle is not less than
$\frac35\pi-\Omega'(\tr)\ge \frac35\pi-\frac{1}{500}$. Let $\tr=\tr ABC$,
the angle $\an B$ being the greatest
one.

As the curvature of the triangle $\tr$ is small, the angles of any lune formed by two
shortests and containing in $\tr$ do not exceed  $\frac{1}{500}$.
Therefore there is a point $X\in AC$ such
that any shortest $BX$ makes angles with $AB$  and $CB$ not less than
$\frac{3}{10}\pi -\frac{1}{500}$.
Some simple calculation shows that the both triangles,  $\tr ABX$ and $\tr CBX$,
satisfy the conditions of Lemma \ref{small-1}. Hence, both of the triangles
are  bi-Lipschitz equivalent to
their comparison triangles respectively. It is easy to see that the same is
true for $\tr ABC$.
Thus  $M$ is bi-Lipschitz equivalent to some surface with a polyhedral metric.

Now it remains to use the fact that if two polyhedral surfaces $X$ and $Y$ are homeomorphic,
then there exists a piecewise linear homeomorphism $g\colon X\rightarrow Y $ between them.
Then according to \cite{RS} we can find such isomorphic subdivisions
%%sub-partitions  ???
$X'$ and $Y'$ of these
triangulations (by flat triangles) that $g$ is linear on every triangle of $g$ and transform
it into corresponding triangle of $Y'$. Taking this into account, we derive that two
homeomorphic polyhedral Alexandrov spaces are bi-Lipschitz equivalent, and
hence the same is valid for two arbitrary homeomorphic Alexandrov spaces satisfying the
conditions of the theorem.

\rightline{$\blacksquare$}

\section {Rotationally symmetric ends}

An end  $T$ is rotationally symmetric if there is an isometry group of
$T$ acting transitively in $\partial T$.
\medskip

{\bf Proof of Theorem \ref{theorem-3}.}

If an end has nonzero speed, then Remark  \ref{Lang-Bonk ends} says
that it is  bi-Lipschitz equivalent to $\R^2$ with a disk removed. We have
mentioned already, that this fact can be proved by a minor modification of the method
of the paper \cite{BonkLang}. Therefore we can suppose that the end  $T$ has zero
speed. In addition, it is sufficient to prove the theorem for polyhedral ends having
a finite number of vertices. And we can restrict ourselves only with ends satisfying
the following conditions:

1) the boundary $\Gamma$ of the end is either a geodesic loop or a polygon, all angles of
which being not greater than $\pi$.

2) $\Omega (T)+\sigma (\Gamma)<\epsilon =\frac{1}{1000}$, where $\sigma (\Gamma )$ is the
variation of the boundary turn of the end.

Indeed,  we can cut off a tubular neighborhood of the boundary in such a way that
the remaining
% reduced ???
end will satisfy the conditions 1) -- 2) and the annulus we cut off will be  bi-Lipschitz
equivalent to a flat annulus (Theorem \ref{theorem-1}).

Further we assume  the conditions 1) -- 2) to be fulfilled. Subsequent proof is
similar to our proof of Lemma \ref{small-tr}. The only difference is that now we will
construct a partition of the end into flat trapezoids; these  trapezoids will be placed
at layers, and all  trapezoids in one layer will have equal highs. Each layer will be
bi-Lipschitz  equivalent to a surface of revolution, bi-Lipschitz constants will be
uniformly bounded, and restrictions of corresponding bi-Lipschitz  maps to boundaries
will preserve lengths.

To do that,  from every ``angular'' point $X_i$, $i=1,\dots ,m,$ of the
geodesic  broken $\Gamma$
(i.e, from points in which  $\Gamma$ has nonzero turn) we emanate a geodesic forming
equal angles with branches of $\Gamma$ starting at $X_i$; i.e., going along the
bisector of the angle between the branches.
Choose  a small number  $h_1>0$ and, in every such a geodesic,   select a point $X_{1i}$
such that
$|X_iX_{1i}|=(\cos\alpha_i)^{-1}h_1$, where  $2\alpha$  is  the turn of
$\Gamma$ at  $X_i$.

Now we connect cyclically the points  $X_{1i}$ by shortests
$X_{1i}X_{1\,i+1}$.
For \linebreak sufficiently small $t_1$ the quadrangles $X_iX_{1i}X_{1\,i+1}X_{i+1}$
are flat nonoverlapping  trapezoids having highs
$h_1$ and angles $\alpha_i$ and $\alpha_{i+1}$ adjoined to the ``bottom'' base.
(If  $\Gamma$ is a loop with zero turn everywhere except  the vertex, we get
one  trapezoid glued  with itself along $X_1X_{11}$.)

Now we will increase  $h_1$ until one of the following events happens:

a) A vertex of the metric appears for the first time in the broken line
$\Gamma_1=X_{11}, X_{12}, \dots ,X_{1m}$.

b) Shortests  $X_iX_{1i}$, $X_{1\,i+1}X_{i+1}$ meet together at a point
$X_{1i}X_{1\,i+1}$ for some $i$.

For sure, several such appearances and meetings can happen simultaneously. It is
not difficult to see that
 nothing but events a) -- b) can happen as curvature is small.

We repeat the described construction starting with the broken line $\Gamma_1$.
(The number
of its vertices may be more numerous just as less numerous than the number of vertices
of the broken line  $\Gamma$.) We get the second row of trapezoids, and so on. Just as
in the proof of Lemma \ref{small-tr}, our supposition guarantees that this process will not
stop until all the vertices of the metric are exhausted. In addition, the turn of the
broken line $\Gamma_k$ is the sum of the turn of the broken line  $\Gamma_{k-1}$
and curvatures of the additional vertices, so that
                   $\sigma(\Gamma_k)\leq \Omega (T)+\sigma(\Gamma)$.
In particular, all the angles of constructed trapezoids differ from  $\pi/2$
not more than $\frac{1}{1000}$.

Denote the length of the broken line $\Gamma_k$, $k=0,1,\dots  $ by  $L_k=L(\Gamma_k)$.
Here we presume  $\Gamma=\Gamma_0$. Now we consider an annular $C_i$ lying on the cone
with the total angle $\theta=t^{-1}_k(L_k-L_{k-1})$ around its vertex $O$ and bounded
by circles of radii $\theta^{-1}L_k$ and $\theta^{-1}L_{k+1}$  centered at $O$.
The lengths of these circles are  $L_k$ and $L_{k+1}$, and they distanced one from another
by $h_k$. ( If $L_k=L_{k+1}$, then the cone degenerates into a cylinder, simplifying
our considerations).

Our next (and final) aim is to show that every annual layer between $\Gamma_k$ and
$\Gamma_{k+1}$ can be bi-Lipschitz mapped onto $C_k$, the constant depending only on
$\epsilon$, we have chosen. The number  $k$ of the layer is supposed to be fixed
and further  it will be omitted in notations. Put  $a_i=|X_{k\,i}X_{k\,i+1}|$,
$b_i=|X_{k+1\,i}X_{k+1\,i+1}|$, $a=\sum a_i$, $b=\sum b_i$. Consider the trapezoid
$ABCD$ with  $|AD|=a$, $|BC|=b$, $\an DAB=\pi$ and mark  points
$A=A_0,A_1,\dots , A_m=D$ and $B=B_0,B_1,\dots ,B_m=C$ on its bases
such that $|A_iA_{i+1}|=a_i$,
$|B_iB_{i+1}|=b_i$. Simple calculations show that angles $\beta_i$ between the
intervals $A_iB_i$ and the  trapezoid bases are uniformly separated from zero
and $\pi$ by the
constant depending only on $\epsilon$.
In fact, as we have
$|a_j-b_j|=h|\tg\alpha_j+\tg\alpha_{j+1}|\leq 2h(|\alpha_j|+|\alpha_{j+1}|),$
(the last inequality holds because  $\epsilon$ is small),
then

\centerline{$|\ctg\beta_i|=
h^{-1}|\sum_{j=1}^{i+1}(a_j-b_j)|\leq 4\sum_{j=1}^m |\alpha_j|\leq\epsilon$.}

Now it is not difficult to construct a map of each trapezoid of
the end partition onto corresponding trapezoid
$A_iB_iA_{i+1}B_{i+1}$ and hence a map of each annual layer onto
corresponding trapezoid $ABCD$. At last, it is easy to map the
trapezoid $ABCD$ onto the corresponding  annual $C_i$ on the cone.
(All the maps can be  constructed in such a way that their
restrictions on the boundaries preserve the lengths).

\rightline{$\blacksquare$}

{\bf Final Remarks: Alexandrov polyhedra.}

Together with Alexandrov surfaces it is possible to consider
 polyhedra glued of them.

 Alexandrov polyhedron is a connected  2-polyhedron  $P$ glued from a finite number of
Alexandrov surfaces under the conditions (i) the gluing are made along the boundary
curves, and (ii) parts attached
one with another have equal lengths (more precisely, the gluing is made along
isometries of boundaries).

A particular case of Alexandrov polyhedra are 2-polyhedra of the curvature bounded
above investigated in \cite{Bur-Buyalo}. Notions of essential edge and maximal face
used there are of pure topological nature and hence can be applied in our case. Each
maximal face is an Alexandrov surface (by Alexandrov's Gluing  Theorem).
We do not exclude a "boundary curve" consisting of one point
(the turn of such a curve is, by definition,
$2\pi-\theta$, where $\theta$ is the total angle around this point). It might be also
possible to allow existence of edges without any faces adjacent to them.
Nevertheless we will not do it here. Also for simplicity  we restrict
ourselves by  polyhedra without boundary edges.
\medskip

{\bf Remark.} This definition of Alexandrov polyhedra is of
constructive \linebreak character.
We would like to find an axiomatic (workable) definition of  Alexandrov polyhedra.
But alas, we could not manage to do it.

\medskip

For  Alexandrov polyhedra, curvature can be naturally defined. In fact, it is sufficient to
define it for essential vertices  and on subsets of essential edges, as
out of these vertices and edges curvature of a set is assumed to be equal to its curvature
as a part of the maximal face. If $\gamma$ is an essential edge,
$e\subset \gamma$, and $\tau_i$ is the turn of  $\gamma$ as of the part of
the boundary of the $i$-th face adjacent to $\gamma$, then by definition
$\omega (e)=\sum_i\tau_i(e)$, where sum is spread over all faces adjacent to the edge
$\gamma$. For the vertex $p$, curvature is defined by the formula
$$
\omega (p)=(2-\chi(\Sigma_p))\pi-s(\Sigma_p),
$$
where $\chi$ and  $s$ are Euler characteristic and the length of the graph $\Sigma_p$.
Under this definition of the curvature, for compact Alexandrov polyhedra Gauss--Bonnet
Theorem is true in its ordinary form.

The situation with positive and negative parts of the curvature is more complicate.
We also need to make the definition of ends more precise.
Having in mind the straightforward
generalizations of the theorem above, we must exclude analogies of peak points
and the ends with zero speed growth. The former means that we should suppose
that the maximal faces having to contain neither inner points with the curvature $2\pi$,
nor boundary points for which the boundary has the turn $\pi$.
This restriction can be described more strictly using the language of curvatures
of the essential edges and vertices. Namely, we need to suppose that there are
no points $p$ in which $\max_i\tau^+_i(p)=\pi$ on the essential edges. Here maximum is
taken over all the faces adjacent  to $\gamma$. For the vertex $p$, let us consider
its link $\Lambda(p)$.

Vertices of a link correspond to essential edges of the polyhedra (starting at $p$)
and its edges correspond to the maximal faces to which these edges are adjacent;
in particular,
the edges of a link can be circles (do not containing essential vertices or containing
only one vertex).
The space of directions at the point $p$ naturally induces a semi-metric on $\Lambda(p)$.
Now let us assume this semi-metric to be a metric; i.e., $\Lambda(p)$ to be
homeomorphic to
the space of directions. In other words, we claim that every edge of a link has nonzero
length in the angle metric.

Now about ends. Let  $Q$ be a part of a polyhedron with complete infinite metric
homeomorphic to the direct product of a finite graph and  $\R^+$.
It is reasonable  to consider
a cone; i.e., single-point compactification $\bar Q=Q\cup \{p\}$ of the space $Q$,
so that all  paths going to the point $p$ have infinite length.

We will think of ends as of subcones of this cone which are cones over the arcs of
the graph connecting its two neighboring essential vertices.
 Thus an end can look either as $S^1\times \R^+$ or as
 $[a,b]\times R^+$.

With these definitions in mind the theorems proved above can be easily extended from
surfaces to the  Alexandrov polyhedra.


\begin{thebibliography}{99}

\bibitem [AZ]{AlexZalg} A. D. Alexandrov and V. A. Zalgaller,
\textit{Intrinsic Geometry of Surfaces},
AMS Transl. Math. Monographs, \textbf{15} (1967), Providence, RI (transl. from Russian).

\bibitem [Bak]{Bakelman}
И. Я. Бакельман,
\textit{Чебышевские сети на многообразиях ограниченной
кривизны}, Труды МИАН, \textbf{76} (1965), стр. 124-129 (in Russian).

\bibitem [BB]{Bur-Buyalo}
Yu. Burago and S. Buyalo, \textit{Metrics of curvature bounded above on
2-polyhedra. II,} St. Petersburg Math. J. \textbf{10}, no. 4, (1998), 62-112.

\bibitem[BL]{BonkLang}
M. Bonk, U. Lang, \textit{Bi-Lipschitz parametrization of
surfaces},  Math. Ann. \textbf{327}, 2003, 135 - 169, (DOI:
10.1007/s00208-003-0443-8).


\bibitem[Gr]{Grieser} D. Grieser,
\textit{Quasiisometry of singular metrics},
Houston J. of Mathem., \textbf{28}, no. 4 (2002), 741-752.

\bibitem[Hub]{Huber} A. Huber,
\textit{On subharmonic functions and differencial geometry in the
large}, Comment. Math. Helv, \textbf{32},  (1957), 13 - 72.

\bibitem [Resh]{Reshetnyak}
Yu. Reshetnyak, {Two dimensional manifolds of bounded curvature}, pp. 3-163.
In: Yu. G. Reshetnyak (Ed.), Geometry IV, Encyclopaedia of Math. Sci.,
\textbf{70}, Springer 1993.

\bibitem [RS]{RS} C. P. Rourke and B. J. Sanderson,
\textit{Introduction to Piecewise-Linear Topology}, Springer, 1972
in ser. ``Ergebnisse der Mathematik und ihrer Grenzgebiete'', Band 69.

\bibitem[Ver]{Verner} А. Вернер,
\textit{Условие конечносвязности полных незамкнутых surfaceей},
Уч. зап. ЛГПИ им. А. И. Герцена, \textbf{395}, no. 4 (1970),
100 - 131 (in Russian).

\end{thebibliography}
\end{document}